\documentclass[a4paper,12pt]{amsart}

\usepackage{latexsym}
\usepackage{amssymb}
\usepackage{amsthm}
\usepackage{amscd}
\usepackage{amsmath}
\usepackage{graphics}
\usepackage[all]{xy}

\newtheorem{Theorem}{Theorem}[section]

\newtheorem{Proposition}[Theorem]{Proposition}

\newtheorem{Definition}[Theorem]{Definition}

\newtheorem{Question}[Theorem]{Question}

\def\fka{{\mathfrak a}}
\def\fkb{{\mathfrak b}}

\def\fkm{{\mathfrak m}}

\def\opn#1#2{\def#1{\operatorname{#2}}}

\opn\Spec{Spec}
\opn\Supp{Supp}
\opn\supp{supp}
\opn\Max{Max}
\opn\max{max}
\opn\Min{Min}
\opn\min{min}
\opn\Ass{Ass}
\opn\Assh{Assh}
\opn\depth{depth}
\opn\rank{rank}

\opn\Mat{Mat}

\opn\Tot{Tot}

\opn\Sym{Sym}
\def\Rees{{\mathcal R}}

\def\S{{\mathcal S}}

\def\ZZ{\mathbb{Z}}  

\opn\div{div}
\opn\Div{Div}
\opn\cl{cl}
\opn\Cl{Cl}

\opn\Ker{Ker}
\opn\Coker{Coker}
\opn\Im{Im}
\opn\Hom{Hom}
\opn\Tor{Tor}
\opn\Ext{Ext}
\opn\End{End}
\opn\Fitt{Fitt}
\opn\Aut{Aut}
\opn\id{id}

\opn\nat{nat}
\opn\pff{pf}
\opn\Pf{Pf}
\opn\GL{GL}
\opn\SL{SL}
\opn\G{G}
\opn\E{E}
\opn\H{H}
\opn\M{M}
\opn\mod{mod}
\opn\ord{ord}
\opn\det{det}
\opn\Soc{Soc}

\opn\chara{char}
\opn\length{\ell}
\opn\pd{pd}
\opn\rk{rk}
\opn\projdim{proj\,dim}
\opn\injdim{inj\,dim}
\opn\rank{rank}
\opn\depth{depth}
\opn\grade{grade}
\opn\height{ht}
\opn\embdim{emb\,dim}
\opn\codim{codim}

\renewcommand{\tilde}{\widetilde}

\renewcommand{\hat}{\widehat}

\title{A Note on the Buchsbaum-Rim multiplicity of a parameter module}

\author{Futoshi Hayasaka}
\author{Eero Hyry}

\address{Department of Mathematics, School of Science and Technology, 
Meiji University, 1-1-1 Higashimita, Tama-ku, Kawasaki 214--8571, JAPAN}
\email{hayasaka@isc.meiji.ac.jp}

\address{Department of Mathematics and Statistics, University of Tampere, 
33014 Tampereen yliopisto, FINLAND}
\email{Eero.Hyry@uta.fi}

\date{\today}


\begin{document}



\begin{abstract}
In this article we prove that the Buchsbaum-Rim multiplicity $e(F/N)$ of a parameter module $N$ in a free module $F=A^r$ is bounded above by the colength $\ell_A(F/N)$. Moreover, we prove that once the equality $\ell_A(F/N)=e(F/N)$ holds true for some parameter module $N$ in $F$, then the base ring $A$ is Cohen-Macaulay. 
\end{abstract}

\maketitle

\section{Introduction}

Let $(A, \fkm)$ be a Noetherian local ring with the maximal ideal $\fkm$ and $d=\dim A>0$. Let $F=A^r$ be a free module of rank $r>0$, and let $M$ be a submodule of $F$ such that $F/M$ has finite length and $M \subseteq \fkm F$. 

In their article \cite{BR2} from 1964 Buchsbaum and Rim introduced and studied a multiplicity associated to a submodule of finite colength in a free module. This multiplicity, which generalizes the notion of Hilbert--Samuel multiplicity for ideals, is nowadays called the Buchsbaum--Rim multiplicity. In more detail, it first turns out that the function $$\lambda(n):=\ell_A(\S_{n}(F)/\Rees_{n}(M))$$ is eventually a polynomial of degree $d+r-1$, where $\S_A(F)=\bigoplus_{n \geq 0}\S_{n}(F)$ is the symmetric algebra of $F$ and $\Rees(M)=\bigoplus_{n \geq 0} \Rees_{n}(M)$ is the image of the natural homomorphism from $\S_A(M)$ to $\S_A(F)$. The polynomial $P(n)$ corresponding to $\lambda(n)$ can then be written in the form 
$$P(n)=\sum_{i=0}^{d+r-1}(-1)^i e_i \binom{n+d+r-2-i}{d+r-1-i} $$ with integer coefficients $e_i$. The {\it{Buchsbaum-Rim multiplicity of}} $M$ {\it{in}} $F$, denoted by $e(F/M)$, is now defined to be the coefficient $e_0$.  

Buchsbaum and Rim also introduced in their article the notion of a parameter module (matrix), which generalizes the notion of a parameter ideal (system of parameters). The module $N$ in $F$ is said to be {\it{a parameter module in}} $F$, if the following three conditions are satisfied: (i) $F/N$ has finite length, (ii) $N \subseteq \fkm F$, and (iii) $\mu_A(N)=d+r-1$, where $\mu_A(N)$ is the minimal number of generators of $N$. 

Buchsbaum and Rim utilized in their study the relationship between the Buchsbaum-Rim multiplicity and the Euler-Poincar\'e characteristic of a certain complex, and proved the following: 

\begin{Theorem}[Buchsbaum-Rim {\cite[Corollary 4.5]{BR2}}]\label{CM case}
Let $(A, \fkm)$ be a Noetherian local ring of dimension $d>0$. Then the following statements are equivalent:
\begin{enumerate}
\item[$(1)$] $A$ is a Cohen-Macaulay local ring; 
\item[$(2)$] For any rank $r>0$, the equality 
$\ell_A(F/N)=e(F/N)$ holds true for every parameter module $N$ in $F=A^r$. 
\end{enumerate}
\end{Theorem}

Then it is natural to ask the following:

\begin{Question}\label{ques} \

\begin{enumerate}
\item[$(1)$] Does the inequality $\ell_A(F/N) \geq e(F/N)$ hold true for any parameter module $N$ in $F$?
\item[$(2)$] Does the equality $\ell_A(F/N)=e(F/N)$ for some parameter module $N$ in $F$ imply that the ring $A$ is Cohen-Macaulay?
\end{enumerate}

\end{Question}

The purpose of this article is to give a complete answer to Question \ref{ques}. Our results can be summarized as follows:

\begin{Theorem}\label{main}
Let $(A, \fkm)$ be a Noetherian local ring of dimension $d>0$. 
\begin{enumerate}
\item[$(1)$] For any rank $r>0$, the following two inequalities 
$$\ell_A(F/N) \geq e(F/N) \ {\mbox{and}} \ \ell_A(A/I(N)) \geq e(F/N) $$ always hold true for every parameter module $N$ in $F=A^r$, where $I(N)$ is the $0$-th Fitting ideal of $F/N$. 
\item[$(2)$] The following statements are equivalent: 
\begin{enumerate}
\item[$(i)$] $A$ is a Cohen-Macaulay local ring; 
\item[$(ii)$] For some rank $r>0$, there exists a parameter module $N$ in $F=A^r$ such that the equality $\ell_A(F/N)=e(F/N)$ holds true; 
\item[$(iii)$] For some rank $r>0$, there exists a parameter module $N$ in $F=A^r$ such that the equality $\ell_A(A/I(N))=e(F/N)$ holds true. 
\end{enumerate}
When this is the case, the equality 
$\ell_A(F/N)=\ell_A(A/I(N))=e(F/N)$ holds true for all parameter modules $N$ in $F=A^r$ of any rank $r>0$. 
\end{enumerate}
\end{Theorem}

Note that the equality $\ell_A(F/N)=\ell_A(A/I(N))$ is known by \cite[2.10]{BV2}.

The proof of our Theorem \ref{main} will be completed in section 4. Section 2 is of preliminary character. In that section we will recall the definition and some basic facts about the generalized Koszul complex. In order to prove Theorem \ref{main}, we will investigate in section 3 the higher Euler-Poincar\'e characteristics of the generalized Koszul complex and show that they are non-negative. Finally, in section 4, we will obtain Theorem \ref{main} as a corollary of a more general result (Theorem \ref{general}). 

\section{Preliminaries}

In this section we will recall the definition and some basic facts about the generalized Koszul complex introduced in \cite{BE1, K} (for more details, see also \cite[Appendix A2.6]{E}). 

Let $A$ be a commutative Noetherian ring, and let $n \geq r >0$ be integers. 
Let $\fka=(a_{ij})$ be an $r \times n$ matrix over $A$, and let $I_r(\fka)$ denote the ideal generated by the maximal minors of $\fka$. 
Let $F$ and $G$ be free modules with bases 
$\{f_1, \dots , f_r \}$ and $\{e_1, \dots , e_n\}$, respectively.  
Let $S$ be the symmetric algebra of $F$, and $S_{\ell}$ the $\ell$-th symmetric power of $F$. Let $\wedge$ be the exterior algebra of $G$, and $\wedge^{\ell}$ the $\ell$-th exterior power of $G$. Associated with the $i$-th row 
$[a_{i1}  \cdots \ a_{in}]$ of $\fka$, there is a differentiation homomorphism 
$\delta_i: \wedge \to \wedge$ given by 
$$\delta_i(f_{j_1} \wedge \dots \wedge f_{j_p})=\sum_{k=1}^p(-1)^{k-1}a_{ij_k}f_{j_1} \wedge \dots \wedge \hat{f_{j_k}} \wedge \dots \wedge f_{j_p}. $$
Let $f_i : S \to S$ and $f_{i}^{-1}:S \to S$ denote the multiplication and division maps by $f_i$, respectively, i.e., 
$$f_i^{-1}(f_1^{\mu_1} \cdots f_i^{\mu_i} 
\cdots f_r^{\mu_r})=\left\{
\begin{array}{ll}
f_1^{\mu_1} \cdots f_i^{\mu_i-1} \cdots f_r^{\mu_r} & (\mu_i > 0) \\
0 & (\mu_i=0). 
\end{array}
\right.
$$
Then the generalized Koszul complex $K_{\bullet}(\fka; t)$ 
associated to a matrix $\fka$ and an integer $t$ is the complex 
$$ K_{\bullet}(\fka; t) : \ \ \cdots \to K_{p+1}(\fka;t) \stackrel{d_{p+1}}{\to} K_p(\fka;t) \stackrel{d_{p}}{\to} K_{p-1}(\fka;t) \to \cdots $$
defined by
$$K_p(\fka; t)=\left\{
	\begin{array}{ll}
		\wedge^{r+p-1} \otimes_A S_{p-t-1} & (p \geq t+1) \\
		\wedge^p \otimes_A S_{t-p} & (p \leq t)  
	\end{array}
	\right.$$
and
$$d_{p+1}=\left\{
	\begin{array}{ll}
		\sum_{j=1}^r \delta_j \otimes f_j^{-1} & \ \ (p > t) \\
		\delta_r \circ \cdots \circ \delta_1 \otimes 1 & \ \ (p = t) \\
		\sum_{j=1}^r \delta_j \otimes f_j & \ \ (p < t). 
	\end{array}
	\right. $$

The generalized Koszul complex $K_{\bullet}(\fka; t)$ is a free complex of $A$-modules. 
We note that it is of length $n-r+1$, when $-1 \leq t \leq n-r+1$. Also recall that $K_{\bullet}(\fka; t)$ coincides with the ordinary Koszul complex for any $t$ in the case $r=1$, whereas
$K_{\bullet}(\fka; 0)$ is the Eagon-Northcott complex and $K_{\bullet}(\fka; 1)$ is the Buchsbaum-Rim complex. Moreover, the generalized Koszul complex has the following important properties (see \cite{K, K3} and \cite[Appendix A2.6]{E}):

\begin{Proposition}\label{fact}
Let $\fka$ be an $r \times n$ matrix over $A$ with $n \geq r >0$. Then 
\begin{enumerate}
\item[$(1)$]\cite[Theorem 1]{K} For any $t, p \in \ZZ$, $I_r(\fka) H_p(K_{\bullet}(\fka; t))=(0)$. 
\item[$(2)$]\cite[Theorem A2.10]{E} If the grade of $I_r(\fka)$ is at least $n-r+1$, then 
$K_{\bullet}(\fka; t)$ is acyclic for all $-1 \leq t \leq n-r+1$. Furthermore, 
if $\fka$ is a generic matrix, then $K_{\bullet}(\fka; t)$ is acyclic for all $t \geq -1$.  
\end{enumerate}
\end{Proposition}

\section{Higher Euler-Poincar\'e characteristics}

In this section we will investigate higher Euler-Poincar\'e characteristics of a generalized Koszul complex. 

Throughout this section, let $(A, \fkm)$ be a Noetherian local ring of dimension $d>0$. Let $F=A^r$ be a free module of rank $r>0$ with a basis $\{ f_1, \dots , f_r \}$. Let $M$ be a submodule of $F$ generated by $c_1, c_2, \dots , c_n$, where $n=\mu_A(M)$ is the minimal number of generators of $M$. 
Writing $ c_j=c_{1j}f_1 + \dots + c_{rj}f_r$ for some $c_{ij} \in A$, we have an $r \times n$ matrix $(c_{ij})$ associated to $M$. We call this matrix 
the matrix of $M$, and denote it by $\tilde{M}$. 
Let $I(M)=\Fitt_0(F/M)$ be the $0$-th Fitting ideal of $F/M$. 
We assume that $F/M$ has finite length and $M \subseteq \fkm F$. 
Then $I(M)$ is an $\fkm$-primary ideal, because $\sqrt{I(M)} = \sqrt{\mathrm{Ann}_A(F/M)}$. Hence each homology module $H_p(K_{\bullet}(\tilde{M}; t))$ has finite length by Proposition \ref{fact}(1). So the Euler-Poincar\'e characteristics of $K_{\bullet}(\tilde{M}; t)$ can be defined as follows:

\begin{Definition}
For any integer $q \geq 0$, we set 
$$\chi_q(K_{\bullet}(\tilde{M}; t)):=\sum_{p \geq q}(-1)^{p-q}\ell_A(H_p(K_{\bullet}(\tilde{M}; t)))$$
and call it the $q$-th partial Euler-Poincar\'e characteristic of $K_{\bullet}(\tilde{M}; t)$. When $q=0$, we simply write $\chi(K_{\bullet}(\tilde{M}; t))$ for $\chi_0(K_{\bullet}(\tilde{M}; t))$, and call it the Euler-Poincar\'e characteristic of $K_{\bullet}(\tilde{M}; t)$. 
\end{Definition}

Buchsbaum and Rim studied in \cite{BR2} the Euler-Poincar\'e characteristic of the Buchsbaum-Rim complex in analogy with
the Euler-Poincar\'e characteristic of the ordinary Koszul complex in the case of usual multiplicities. In 1985 Kirby 
investigated in \cite{K2} Euler-Poincar\'e characteristics of the complex $K_{\bullet}(\tilde{M}; t)$ for all $t$ and proved the following: 

\begin{Theorem}[Buchsbaum-Rim, Kirby]\label{chi}
For any integer $t \in \ZZ$, we have 
$$
\chi(K_{\bullet}(\tilde{M}; t))=
\left\{ 
\begin{array}{cl}
e(F/M) & (n=d+r-1), \\
0      & (n>d+r-1), 
\end{array}
\right.
$$
where $n=\mu_A(M)$ is the minimal number of generators of $M$. In particular, $\chi(K_{\bullet}(\tilde{M}; t)) \geq 0$ for all $t \in \ZZ$. 
\end{Theorem}

The last statement holds for the higher Euler-Poincar\'e characteristics, too: 

\begin{Theorem}\label{nonnegative}
For any $q \geq 0$ and any $t \geq -1$, we have 
$$\chi_q(K_{\bullet}(\tilde{M}; t)) \geq 0. $$
\end{Theorem}

\begin{proof}
We use ideas from \cite{BR3}. Let $\tilde{M}=(c_{ij}) \in \Mat_{r \times n}(A)$ be the matrix of $M$ and $X=(X_{ij})$ be the generic matrix of the same size $r \times n$. Let $A[X]=A[X_{ij} \mid 1 \leq i \leq r, \ 1 \leq j \leq n]$ be a polynomial ring over $A$ and let $B=A[X]_{(\fkm, X)}$. We will consider the ring 
$A$ as a $B$-algebra via the substitution homomorphism
$\phi: B \to A \ ; X_{ij} \mapsto c_{ij}$. 
Let $$\fkb=\Ker \phi=(X_{ij}-c_{ij} \mid 1 \leq i \leq r, \ 1 \leq j \leq n)B.$$ 
We note here that $K_{\bullet}(X; t) \otimes_B A \cong K_{\bullet}(\tilde{M}; t)$, because the generalized Koszul complex is compatible with the base change. 
Let $C_t(X):=H_0(K_{\bullet}(X; t))$. By Proposition \ref{fact}(2), the complex $K_{\bullet}(X; t)$ is a $B$-free resolution of the $B$-module $C_t(X)$ for any $t \geq -1$. By tensoring with $A$ and taking the homology, we have that 
\begin{eqnarray*}
H_p(K_{\bullet}(\tilde{M}; t)) & \cong & H_p(K_{\bullet}(X; t) \otimes_B A) \\
                               & \cong & \Tor_p^B(C_t(X), A)
\end{eqnarray*}
for all $p \geq 0$. 
On the other hand, since the ideal $\fkb$ in $B$ is generated by a regular sequence of length $rn$, the ordinary Koszul complex $K_{\bullet}(\fkb)$ associated to the sequence $\tilde{\fkb}$ is a $B$-free resolution of $A$. Hence, by tensoring with $C_t(X)$, we can compute the $\Tor$ as follows:
\begin{eqnarray*}
\Tor_p^B(C_t(X), A) \cong H_p(K_{\bullet}(\fkb) \otimes_B C_t(X)). 
\end{eqnarray*}
Therefore, for any $p \geq 0$, we have 
$$H_p(K_{\bullet}(\tilde{M}; t)) \cong H_p(K_{\bullet}(\fkb) \otimes_B C_t(X)). $$
It follows that for any $t \geq -1$ and any $q \geq 0$ we have the equality 
$$\chi_q(K_{\bullet}(\tilde{M}; t))=\chi_q(K_{\bullet}(\fkb) \otimes_B C_t(X)).$$ Here the right hand side is non--negative by Serre's Theorem (\cite[Ch. IV Appendix II]{Serre}). Therefore $\chi_q(K_{\bullet}(\tilde{M}; t)) \geq 0$. 
\end{proof}

\section{Proof of Theorem \ref{main}}

Theorem \ref{main} will be a consequence of the following more general result:

\begin{Theorem}\label{general}
Let $(A, \fkm)$ be a Noetherian local ring of dimension $d>0$. 
\begin{enumerate}
\item[$(1)$] For any rank $r>0$, the inequality 
$\ell_A(H_0(K_{\bullet}(\tilde{N}; t))) \geq e(F/N)$
holds true for any integer $t \geq -1$ and any parameter module $N$ in $F=A^r$. 
\item[$(2)$] The following statements are equivalent: 
\begin{enumerate}
\item[$(i)$] $A$ is a Cohen-Macaulay local ring; 
\item[$(ii)$] For some rank $r>0$, there exists an integer $-1 \leq t \leq d$ and a parameter module $N$ in $F=A^r$ such that the equality $\ell_A(H_0(K_{\bullet}(\tilde{N}; t)))=e(F/N)$ holds true. 
\end{enumerate}
When this is the case, the equality 
$\ell_A(H_0(K_{\bullet}(\tilde{N}; t)))=e(F/N)$ holds true 
for any integer $-1 \leq t \leq d$ and any parameter module $N$ in $F=A^r$ of any rank $r>0$. 
\end{enumerate}
\end{Theorem}

\begin{proof}

(1): Let $N$ be a parameter module in $F=A^r$, and let $t \geq -1$. By Theorem \ref{chi} we obtain that 
\begin{eqnarray*}
e(F/N) & = & \chi(K_{\bullet}(\tilde{N}; t)) \\
		     & = & \ell_A(H_0(K_{\bullet}(\tilde{N}; t)))-\chi_1(K_{\bullet}(\tilde{N}; t)). 
\end{eqnarray*}
Since $\chi_1(K_{\bullet}(\tilde{N}; t)) \geq 0$ by Theorem \ref{nonnegative}, the desired inequality follows. 

(2): Assume that $A$ is Cohen-Macaulay. Let $N$ be any parameter module in $F=A^r$ of any rank $r>0$. Let $n=\mu_A(N)=d+r-1$. Then $\grade I(N)= \height I(N) = d = n-r+1$. Hence, by Proposition \ref{fact}(2), $K_{\bullet}(\tilde{N}; t)$ is acyclic for all $-1 \leq t \leq n-r+1=d$. Therefore, by Theorem \ref{chi}, we have $e(F/N)=\chi(K_{\bullet}(\tilde{N}; t))=\ell_A(H_0(K_{\bullet}(\tilde{N}; t)))$. This proves the implication (i) $\Rightarrow$ (ii), and also the last assertion.  

It remains to show the implication (ii) $\Rightarrow$ (i). Assume that there exist integers $r > 0$, $-1 \leq t \leq d$, and a parameter module $N$ in $F=A^r$ such that $\ell_A(H_0(K_{\bullet}(\tilde{N}; t)))=e(F/N)$. Arguing as in the proof of Theorem \ref{nonnegative} and using the same notation, we get
\begin{eqnarray*}
\chi_1(K_{\bullet}(\fkb) \otimes_B C_t(X)) &=& 
\chi_1(K_{\bullet}(\tilde{N}; t)) \\ 
& = & \ell_A(H_0(K_{\bullet}(\tilde{N}; t)))-e(F/N) \\               
& = & 0. 
\end{eqnarray*}
We observe here that $\sqrt{\mathrm{Ann}_B C_t(X)}=\sqrt{I_r(X)}$ (see \cite[Lemma 2.7]{R}). Thus $\dim_B C_t(X)=\dim B/I_r(X)=d+(n+1)(r-1)=rn$ (see \cite[(5.12) Corollary]{BV1}). Therefore $\fkb$ is a parameter ideal of $C_t(X)$. Hence we have the equality 
$$\ell_B(C_t(X)/\fkb C_t(X))-e(\fkb; C_t(X)) = \chi_1(K_{\bullet}(\fkb) \otimes_B C_t(X)) = 0, $$
where $e(\fkb; C_t(X))$ is the multiplicity of the module $C_t(X)$ with respect to an ideal $\fkb$. 
Since $\ell_B(C_t(X)/\fkb C_t(X))=e(\fkb; C_t(X))$, this implies that $C_t(X)$ is a Cohen-Macaulay $B$-module. 
On the other hand, $\pd_B C_t(X)=d$, because the complex $K_{\bullet}(X; t)$ is a minimal $B$-free resolution of $C_t(X)$ of length $n-r+1=d$. Hence, by the Auslander-Buchsbaum formula, we have 
\begin{eqnarray*}
d+rn &=& \pd_B C_t(X)+\depth_B C_t(X) \\
     &=& \depth B \\
     & \leq & \dim B \\
     &=& d+rn. 
\end{eqnarray*}
Thus $\depth B=\dim B$ so that $B$ is Cohen-Macaulay. Therefore $A$ is also a Cohen-Macaulay local ring. 
\end{proof}

Taking $t=0, 1$ in Theorem \ref{general}, now readily gives Theorem \ref{main}. 
\bigskip

We want to close this article with a question. For that, let us first recall the notion of Buchsbaum local ring, which was introduced by St\"{u}ckrad and Vogel (for more details on Buchsbaum rings, we refer the reader to \cite{SV}). Let $A$ be a Noetherian local ring. Then $A$ is said to be a {\it{Buchsbaum}} {\it{ring}}, if the difference $$\ell_A(A/Q) - e(A/Q)$$ between the colength and multiplicity of  a parameter ideal $Q$ in $A$ is independent of the choice of $Q$. This difference, which is an invariant of a Buchsbaum ring $A$, is denoted by $I(A)$. The ring $A$ is Cohen-Macaulay if and only if it is Buchsbaum and $I(A)=0$. In this sense, the notion of Buchsbaum ring is a natural generalization of that of Cohen-Macaulay ring. In Theorem \ref{main}, the inequality $\ell_A(F/N) \geq e(F/N)$ for any parameter module $N$ in $F$, is an analogue of the well-known inequality $\ell_A(A/Q) \geq e(A/Q)$ for any parameter ideal $Q$ in $A$. Also, the 
characterization of the Cohen-Macaulay property of $A$ based on the equality $\ell_A(F/N) = e(F/N)$ generalizes the usual one using parameter ideals. With these remarks in mind, it is natural to ask the following question: 

\begin{Question}
Let $F$ be a fixed free module of rank $r>0$. Is it then true that the difference 
$$\ell_A(F/N)-e(F/N)$$ 
between the colength and multiplicity of  a parameter module $N$ in $F$ is independent of the choice of $N$, if the ring $A$ is Buchsbaum?
\end{Question}




\begin{thebibliography}{99}


\bibitem{BV2}
W. Bruns and U. Vetter, 
Length formulas for the local cohomology of exterior powers,
Math. Z. 191 (1986), 145--158

\bibitem{BV1}
W. Bruns and U. Vetter, Determinantal Rings, Lecture Notes in Math.
1327, Springer-Verlag Berlin Heidelberg, 1988

\bibitem{BE1}
D. A. Buchsbaum and D. Eisenbud, Generic free resolutions and a family of generically perfect ideals, Adv. in Math. 18 (1975), 245--301

\bibitem{BR1}
D. A. Buchsbaum and D. S. Rim, 
A generalized Koszul complex, 
Bull. Amer. Math. Soc. 69 (1963), 382--385

\bibitem{BR2}
D. A. Buchsbaum and D. S. Rim, 
A generalized Koszul complex. II. Depth and multiplicity, 
Trans. Amer. Math. Soc. 111 (1964), 197--224

\bibitem{BR3}
D. A. Buchsbaum and D. S. Rim, 
A generalized Koszul complex. III. A Remark on Generic Acyclicity, 
Proc. Amer. Math. Soc. 16 (1965), 555--558


\bibitem{E}
D. Eisenbud, 
Commutative algebra. With a view toward algebraic geometry, 
Graduate Texts in Mathematics, 150, Springer-Verlag, New York, 1995 


\bibitem{K}
D. Kirby, 
A sequence of complexes associated with a matrix, 
J. London Math. Soc. 7 (1974), 523--530

\bibitem{K2}
D. Kirby, 
On the Buchsbaum-Rim multiplicity associated with a matrix, 
J. London Math. Soc. (2) 32 (1985), no. 1, 57--61 

\bibitem{K3}
D. Kirby, Generalized Koszul complexes and the extension functor, 
Comm. Algebra 18 (1990), no. 4, 1229--1244



\bibitem{R}
A. G. Rodicio, 
On the rigidity of the generalized Koszul complexes with applications to Hochschild homology, J. Algebra 167 (1994), no. 2, 343--347

\bibitem{Serre}
J-P. Serre, Local Algebra (Translated from the French by CheeWhye Chin), Springer Monographs in Mathematics, Springer-Verlag Berlin Heidelberg 2000

\bibitem{SV}
J. St\"uckrad and W. Vogel, Buchsbaum rings and applications, Springer-Verlag Berlin Heidelberg New York London Paris Tokyo 1986

\end{thebibliography}
\end{document}